%% file: main.tex
\documentclass[openany,reqno,a4paper,12pt]{amsart}
\input{header.tex}

\usepackage[backend=bibtex,style=numeric]{biblatex}
\bibliography{main.bib}

\usepackage{geometry}
\newcommand{\conn}{\leftrightarrow}
\newcommand{\disconn}{\nleftrightarrow}
\DeclareMathOperator{\supp}{supp}

\begin{document}

\title{A new proof of the bunkbed conjecture in the $p\uparrow 1$ limit}

\author{Lawrence Hollom}
\address{Department of Pure Mathematics and Mathematical Statistics (DPMMS), University of Cambridge, Wilberforce Road, Cambridge, CB3 0WA, United Kingdom}
\email{lh569@cam.ac.uk}



\begin{abstract}
  For a finite simple graph $G$, the bunkbed graph $G^\pm$ is defined to be the product graph $G\square K_2$.
  We will label the two copies of a vertex $v\in V(G)$ as $v_-$ and $v_+$.
  The bunkbed conjecture, posed by Kasteleyn, states that for independent bond percolation on $G^\pm$, percolation from $u_-$ to $v_-$ is at least as likely as percolation from $u_-$ to $v_+$, for any $u,v\in V(G)$.
  Despite the plausibility of this conjecture, so far the problem in full generality remains open. 
  Recently, Hutchcroft, Nizi\'{c}-Nikolac, and Kent gave a proof of the conjecture in the $p\uparrow 1$ limit. 
  Here we present a new proof of the bunkbed conjecture in this limit, working in the more general setting of allowing different probabilities on different edges of $G^\pm$.
\end{abstract}

\maketitle

\section{Introduction}
\label{sec:intro}
In this introduction we state two forms of the bunkbed conjecture, discuss briefly the known results about special cases of the conjecture, and then proceed to state the result that we will prove.

If $G=(V,E)$ is a finite simple graph, then the \emph{bunkbed graph} $G^\pm$ is the box product $G\square K_2$.
Thus, $G^\pm=(V^\pm,E^\pm)$ consists of two copies of $G$, labelled $G^+=(V^+, E^+)$ and $G^-=(V^-,E^-)$, where $V^+\defined\set{v_+\st v\in V}$, and similarly for $V^-$, together with all `vertical' edges connecting $v_+$ and $v_-$ added for all vertices $v\in V$. That is,
\begin{align*}
    E(G^\pm)=E^-\union E^+\union\set{v_-v_+\st v\in V}.
\end{align*}

For a graph $G$, we call a function $\bfp\st E(G)\rightarrow [0,1]$ an \emph{edge-probability function}.
We call an edge-probability function on a bunkbed graph $G^\pm$ \emph{symmetric} if for every edge $xy\in E$ we have $\bfp(x_+y_+)=\bfp(x_-y_-)$.

Given an edge-probability function $\bfp$ for a graph $G$, we can define a percolation process as follows.
We set each edge $e\in E(G)$ to be \emph{open} with probability $\bfp(e)$ and \emph{closed} with probability $1-\bfp(e)$, independently of all other events.
We then define $\Prob_\bfp(u\conn v)$ to be the probability that there is a path of open edges between two vertices $u$ and $v$ of $G$.
As usual we call this probability the two-point function.

As discussed by Rudzinski \cite{Rudzinski2016} there are several different statements referred to as `the bunkbed conjecture'. We first state the conjecture in its most general form.

\begin{conjecture}[Bunkbed conjecture]
\label{conj:bunkbed-general}
Let $G=(V,E)$ be a finite simple graph with bunkbed graph $G^\pm$, and let $\bfp\st E(G^\pm)\rightarrow [0,1]$ be a symmetric edge-probability function. Then for any $u,v\in V$, we have
\begin{align}
    \Prob_\bfp(u_-\conn v_-)\geq \Prob_\bfp(u_-\conn v_+).
\end{align}
\end{conjecture}

The special case in which $\bfp$ takes some constant value $p\in [0,1]$ is \emph{Bernoulli bond percolation}.
In this case we denote the two point function by $\Prob_p(u\conn v)$, and have the following weaker form of the conjecture.

\begin{conjecture}[Uniform bunkbed conjecture]
\label{conj:bunkbed-uniform}
Let $G=(V,E)$ be a finite simple graph with bunkbed graph $G^\pm$. Then for any $u,v\in V$ and any $p\in[0,1]$, we have
\begin{align}
    \Prob_p(u_-\conn v_-)\geq \Prob_p(u_-\conn v_+).
\end{align}
\end{conjecture}

The bunkbed conjecture was first posed by Kasteleyn in 1985, in the form of Conjecture \ref{conj:bunkbed-general}, as noted by Kahn and Berg \cite[Remark 5]{KahnBerg2001}. 
Despite attracting significant interest since then, it so far remains unproved.

Several special cases are however known; de Buyer \cite{Buyer2016} proved Conjecture \ref{conj:bunkbed-uniform} when $G$ is a complete graph and $p=1/2$, and then extended this result to all $p\geq 1/2$ \cite{Buyer2018}.
A little later van Hintum and Lammers \cite{VanHintum2019} used a different method to resolve the case of the complete graph for all $p\in [0,1]$.

Further developments have been made by Richthammer \cite{Richthammer2022}, who proved Conjecture \ref{conj:bunkbed-general} for several other classes of graph, including complete bipartite graphs, complete graphs minus a complete subgraph, and symmetric complete $k$-partite graphs.

More recently, the case of the $p\uparrow 1$ limit has been explored by Hutchcroft, Nizi\'{c}-Nikolac, and Kent \cite{Hutchcroft2021}.
Their results build on the work of Linusson \cite{linusson2011, linusson2019}, and their paper provides a thorough and detailed account of a proof of the following theorem.

\begin{theorem}
\label{thm:bunkbed-uniform-limit}
    For any graph $G$ there is a constant $\eps_0(G)>0$ such that the uniform bunkbed conjecture holds for all $p\geq 1-\eps_0(G)$.
\end{theorem}

Our goal here is to generalise Theorem \ref{thm:bunkbed-uniform-limit} to the setting of Conjecture \ref{conj:bunkbed-general}. This generalisation is stated as follows.

\begin{theorem}
\label{thm:bunkbed-general-limit}
    For any graph $G$ there is a constant $\eps_1(G)>0$ such that the bunkbed conjecture holds for all symmetric edge-probability functions $\bfp$ satisfying $\bfp(e)\geq 1-\eps_1(G)$ for all $e\in E(G^\pm)$.
\end{theorem}

Our proof of Theorem \ref{thm:bunkbed-general-limit} will be elementary and entirely self-contained, proceeding by way of considering cuts in the graph $G^\pm$.


\section{Definitions and Terminology}
\label{sec:definitions}

We call a set $A\sseq V$ a \emph{$(u-v)$-cut} if $A$ contains precisely one of $u$ and $v$. 
Define the \emph{edge-boundary} of $A$ to be $\partial A\defined\set{xy\in E\st x \in A, y\notin A}$.
For $x\in V$ and cut $A\sseq V^\pm$, we define

\begin{align}
\label{def:height}
    h_x(A)\defined \abs{\set{x_-,x_+}\inter A}.
\end{align}

In many cases we will care only about the value of $h_x$ modulo 2.
With this in mind, we define the \emph{support} of a cut as follows.

\begin{align}
\label{def:support}
    \supp(A)\defined\set{x\in V\st h_x(A)=1}.
\end{align}

We define the support of a collection of cuts to be the union of their supports:

\begin{align}
    \supp(A_1,\dotsc,A_r)\defined\supp(A_1)\union\cdots\union\supp(A_r).
\end{align}

Note that a cut and its complement have the same edge-boundary, and so we may assume that a $(u,v)$-cut $A$ contains $u$ and not $v$.
Let $\rS^-$ and $\rS^+$ be the sets of $(u_--v_-)$-cuts containing $u_-$ and $(u_--v_+)$-cuts containing $u_-$ respectively, and let $\rS^\pm\defined\rS^-\union\rS^+$.

We will often consider the connectedness of the support of $A$, and so define
\begin{align}
    \rT^\pm\defined\set{A\in \rS^\pm\st u,v\text{ are in the same component of }G[\supp(A)]}.
\end{align}
with $\rT^-$ and $\rT^+$ defined similarly in terms of $\rS^-$ and $\rS^+$ respectively.

To relate cuts to percolation, note that two vertices $u,v\in V$ are disconnected after percolation if and only if there is a $(u-v)$-cut $A$ such that all edges in $\partial A$ are closed. 
An example of such an $A$ is precisely those vertices reachable from $u$ by paths of open edges.
For a cut $A$, let $E_A$ be the event that all edges in $\partial A$ are closed.

We will be considering our expressions as polynomials, and will say that polynomial $g$ \emph{strictly divides} $h$ if $g$ divides $h$ and $g\neq h$. 


\section{Proof of Theorem \ref{thm:bunkbed-general-limit}}
\label{sec:proof}
After some preliminary work, we reduce proving Theorem \ref{thm:bunkbed-general-limit} to three claims, which we prove in turn. 
First we may rewrite the probability that there is no percolation in terms of cuts, considering first the case of $u_-$ and $v_-$ to find

\begin{align}
\label{eqn:union-expansion}
    \Prob_\bfp(u_-\disconn v_-)=\Prob_\bfp\Bigl(\bigcup _{A\in \rS^-}E_A\Bigr).
\end{align}

We can expand the union into an alternating sum of intersections via the well-known inclusion-exclusion principle.

\begin{align}
\label{eqn:pie}
    \Prob_\bfp\Bigl(\bigcup _{A\in \rS^-}E_A\Bigr)=\sum_{A\in \rS^-}\Prob_\bfp(E_A)-\sum_{A_1\neq A_2\in \rS^-}\Prob_\bfp(E_{A_1}\inter E_{A_2})+\cdots,
\end{align}
with analogous results for $\rS^+$. 
Note that, setting $\bfq(e)\defined 1-\bfp(e)$ for all $e\in E(G^\pm)$, we have
\begin{align}
    \Prob_\bfp(E_{A_1}\inter\cdots\inter E_{A_r})=\bfq(e_1)\cdots\bfq(e_k),
\end{align}
where $\partial A_1\union\cdots\union \partial A_r=\set{e_1,\dotsc,e_k}$.
This gives us an immediate method for interpreting the expressions like those on the right-hand side of equation \eqref{eqn:pie} as polynomials with variables $\bfq(e)$.

To prove our result, we need to show for fixed $u$ and $v$ that $\Prob_\bfp(u_-\disconn v_+)-\Prob_\bfp(u_-\disconn v_-)$ is positive for vertices $u\ne v$.
To this end we now fix $u$ and $v$ and define
\begin{align}
\label{eqn:target}
    f(\bfq)\defined \Prob_\bfp(u_-\disconn v_+)-\Prob_\bfp(u_-\disconn v_-),
\end{align}
where $f(\bfq)$ is considered as a polynomial with variables $\set{\bfq(e)\st e\in E(G^\pm)}$.

Noting that the limit $\bfp\uparrow 1$ is exactly the same as $\bfq \downarrow 0$, we wish to show that every monomial in $f(\bfq)$ with negative coefficient is strictly divided by some monomial in $f(\bfq)$ with positive coefficient.

We are now ready to state the three claims which will prove Theorem \ref{thm:bunkbed-general-limit}. The first two claims deal with terms arising from the first summation of equation \eqref{eqn:union-expansion}.

\begin{claim}
\label{claim:cancellations}
    There is a bijection $\phi\from \rS^+\setminus\rT^+\to \rS^-\setminus\rT^-$ which preserves the value of $\Prob_\bfp(E_A)$. Thus
    \begin{align}
    \label{eqn:cancellation}
        \sum_{A\in{\rS^+\setminus\rT^+}}\Prob_\bfp(E_A)-\sum_{A\in{\rS^-\setminus\rT^-}}\Prob_\bfp(E_A)=0.
    \end{align}
\end{claim}

\begin{claim}
\label{claim:first-order}
    For every cut $A\in\rT^-$ there is a $B\in\rT^+$ such that the polynomial $\Prob_\bfp(E_B)$ strictly divides $\Prob_\bfp(E_A)$.
\end{claim}

The third claim deals with the terms arising from the summations in equation \eqref{eqn:union-expansion} involving more than one cut.

\begin{claim}
\label{claim:higher-order}
    For any integer $r\geq 2$, each nonzero monomial in the polynomial
    \begin{align}
    \label{eqn:higher-order}
        \sum_{\text{distinct }A_1,\dotsc,A_r\in\rS^+}\Prob_\bfp(E_{A_1}\inter\cdots\inter E_{A_r})-\sum_{\text{distinct }A_1,\dotsc,A_r\in\rS^-}\Prob_\bfp(E_{A_1}\inter\cdots\inter E_{A_r})
    \end{align}
    is strictly divided by $\Prob_\bfp(E_B)$ for some cut $B\in\rT^+$.
\end{claim}

Note that once the claims are proved, Theorem \ref{thm:bunkbed-general-limit} follows, as we can rewrite equation \eqref{eqn:target} using equations \eqref{eqn:union-expansion} and \eqref{eqn:pie}, then apply the results of the claims to write
\begin{align*}
    f(\bfq)=\sum_{C\in\rT^+}\Prob_\bfp(E_C)(1+g_C(\bfq)),
\end{align*}
for some polynomials $g_C(\bfq)$, each with constant coefficient zero.
Then as $\bfq\downarrow 0$, the dominating term in $1+g_C(\bfq)$ is $1$, and so $f(\bfq)>0$ for $\bfq$ sufficiently small.

We now prove our three claims.


\subsection{Proof of Claim \ref{claim:cancellations}}
\label{subsec:claim-cancellations}
Take some cut $A\in\rS^+\setminus\rT^+$.
Firstly, if $v\notin\supp(A)$, then we must have $h_v(A)=2$, so $A\in\rS^\pm\setminus\rT^\pm$ and we can set $\phi(A)=A$.

Otherwise, as $G[\supp(A)]$ has $u$ and $v$ in different components, we can map $A$ to the cut $B\in\rS^-\setminus\rT^-$ by swapping $x_+$ with $x_-$ for every $x$ in the $G[\supp(A)]$-component containing $v$.
Note that $\supp(A)=\supp(B)$, so this map is invertible.

This gives us a bijection between $\rS^+\setminus\rT^+$ and $\rS^-\setminus\rT^-$ which preserves the value of $\Prob_\bfp(E_A)$, as required.


\subsection{Proof of Claim \ref{claim:first-order}}
\label{subsec:claim-first-order}
Take any $A\in\rT^-$, and consider the cut $B$ defined as follows.
\begin{align}
    B\defined\set{x_-\st x\in\supp(A)}.
\end{align}

First, we may note that as $u$ and $v$ are in the same component of $G[\supp(A)]$, $B\in\rT^+$.
It remains to show that the polynomial $\Prob_\bfp(E_B)$ strictly divides $\Prob_\bfp(E_A)$.

Take an edge $e\in\partial B$. 
If $e=x_-x_+$ for some $x\in V$, then $x_-\in B$, so $x\in\supp(A)$ and so $e\in\partial A$ as well.

If $e=x_-y_-$ for some $x,y\in V$ with $x\in\supp(A)$ and $y\notin\supp(A)$, then either $x_-y_-$ or $x_+y_+$ is in $\partial A$, as required.

Note that no edge of the form $x_+y_+$ is in $\partial B$, so we have covered all cases and shown that $\Prob_\bfp(E_B)$ divides $\Prob_\bfp(E_A)$.
It remains to see that this division is strict.

Take a path $u=w^{(1)},w^{(2)},\dotsc,w^{(l)}=v$ from $u$ to $v$ in $\supp(A)$, so for each $i$ exactly one of $w^{(i)}_-$ and $w^{(i)}_+$ is present in $A$.
Then if $q_i\defined\bfq(w^{(i)}_-w^{(i+1)}_-)=\bfq(w^{(i)}_+w^{(i+1)}_+)$ are the variables corresponding to the horizontal edges in this path, then we know that no $q_i$ divides the polynomial $\Prob_\bfp(E_B)$.
However, $u_-\in A$ and $v_+\in A$, so for some $j$ we must have $w^{(j)}_-\in A$ and $w^{(j+1)}_+\in A$.
But then $q_j^2$ divides $\Prob_\bfp(E_A)$, so the division is indeed strict, as required.


\subsection{Proof of Claim \ref{claim:higher-order}}
\label{subsec:claim-higher-order}
We may begin as in the proof of Claim \ref{claim:cancellations}.

Let $S\defined\supp(A_1,\dotsc,A_r)$.
If $u$ and $v$ are in different components of $G[S]$, then we can flip the sign of the component of $v$ in all of $A_1,\dotsc A_r$, and the corresponding terms cancel from the sums in equation \eqref{eqn:higher-order}.
Thus we are left to consider collections of cuts with $u$ and $v$ in the same $G[S]$-component.
Now define a cut $B$ as follows.
\begin{align}
    B\defined \set{x_-\st x\in S}.
\end{align}
As in the proof of Claim \ref{claim:first-order}, we know that $B\in\rT^+$.
Note that if $x_-x_+\in\partial B$ is a vertical edge, then $x_-\in B$ and so $x\in S$.
Thus $h_x(A_i)=1$ for some $i$, and so $x_-x_+\in\partial A_i$.

If $x_-y_-\in\partial B$, then assume $x_-\in B$ and $y_-\notin B$.
Then for some $i$ we have $h_x(A_i)=1$, and also $h_y(A_i)\neq 1$ (as the latter holds for all $i$).
Thus either $x_-y_-\in\partial A_i$ or $x_+y_+\in\partial A_i$.
Therefore $\Prob_\bfp(E_B)$ divides $\Prob_\bfp(E_{A_1}\inter\cdots\inter E_{A_r})$.

So it suffices to prove that either the polynomial $\Prob_\bfp(E_{A_1}\inter\cdots\inter E_{A_r})$ strictly divides $\Prob_\bfp(E_B)$, or there is a further cut $C\in\rT^+$ with $\Prob_\bfp(E_C)$ strictly dividing $\Prob_\bfp(E_B)$.

If $\supp(B)=S$ is not a connected subset of $G$, then we can define $C$ to be $B$ restricted to the component of $S$ containing $u$ and $v$.
Then $C\in\rT^+$ and $\Prob_\bfp(E_C)$ strictly divides $\Prob_\bfp(E_B)$, as required.

So now assume that $S$ is connected, and assume further for contradiction that $\Prob_\bfp(E_{A_1}\inter\cdots\inter E_{A_r})=\Prob_\bfp(E_B)$.

First observe that each $A_i$ must be constant on $S$ and constant on $V\setminus S$, as otherwise we would have a horizontal edge in $\partial A_i$ with no corresponding edge in $\partial B$.
Next, note that as $u_-\in A_i$ and $u_+\notin A_i$, all $A_i$ are equal to $B$ on the support of $B$.
Thus the only remaining way in which the $A_i$ can differ is by taking different, constant values on $V\setminus S$.

Assume that $x_-y_-\in\partial B$, with $x\in S$ and $y\notin S$, and consider $h_y(A_i)$.

If $h_y(A_i)=1$, then $y_-y_+\in\partial A_i$, but this vertical edge is not in $\partial B$, contradicting our assumption of equality.
Thus, as $r\geq 2$, the only remaining case is that $r=2$, and wlog $h_x(A_1)=0$ and $h_x(A_2)=2$.
But then $x_-y_-\in\partial A_1$ and $x_+y_+\in\partial A_2$, giving $\bfq(x_-y_-)^2\divides \Prob_\bfp(E_{A_1}\inter E_{A_2})$, again contradicting our assumption of equality.

Thus for any distinct $A_1,\dotsc,A_r$, there is a $C\in\rT^+$ for which the polynomial $\Prob_\bfp(E_{A_1}\inter\cdots\inter E_{A_r})$ strictly divides $\Prob_\bfp(E_C)$.
This is exactly the result we require, so we deduce Claim \ref{claim:higher-order}, and hence also Theorem \ref{thm:bunkbed-general-limit}.


\section{Concluding Remarks}
\label{sec:conclusion}
For the sake of simplicity, we have not bounded the largest value $\eps_1(G)$ for which we can prove Theorem \ref{thm:bunkbed-general-limit}.
To find such a bound using the methods presented here, one would need to count, for each $B\in\rT^+$, the number of terms associated with $B$ by Claims \ref{claim:first-order} and \ref{claim:higher-order}.
The number of terms in the sums in equation \eqref{eqn:higher-order} is doubly exponential in $n$, and the bound attained by our proofs as presented would reflect this.
However, the methods detailed here could be extended to investigate how terms cancel between the sums in \eqref{eqn:higher-order}, which could lead to a much better bound on $\eps_1(G)$. 

A result of Rudzinski and Smyth \cite{Rudzinski2016} shows that if a universal bound $\eps_1(G)>\delta>0$ could be established, then Conjecture \ref{conj:bunkbed-uniform} would be proved for all $p$. 
However, achieving such a bound seems, at present, extremely difficult.

As detailed in \cite[Section 1.1]{Hutchcroft2021}, if we take our graph to be $G=P_n$, the path on $n+1$ vertices, set all $\bfq(e)=q$, and take $u$ and $v$ to be opposite ends of the path, then our polynomial $f(\bfq)$ as in equation \eqref{eqn:target} is equal to $q^{n+1}(1-q)^n$.
In the case of $n=2$, we have $f(q)=q^5-2q^4+q^3$, which cannot be proved to be positive by pairing terms off with each other. 
To extend our method to all values of $\bfq$, we would need a technique to collect and factorise sets of terms, so we could write $f(\bfq)$ as a sum of positive factorised expressions, as in the example above.


\section{Acknowledgement}
\label{sec:acknowledgement}

The author would like to thank Professor B\'{e}la Bollob\'{a}s for his thorough reading of the manuscript and many valuable comments.


\printbibliography


\end{document}

%% file: header.tex
\usepackage{amsmath,amssymb,amsthm}

\usepackage{mathrsfs}

\usepackage{mathtools}
\usepackage{commath}

\usepackage{thmtools}
\usepackage{thm-restate}

\usepackage{cases}
\usepackage{enumitem}
\setlist[enumerate,1]{label=(\roman*)}

\usepackage[pdftex, pdfborderstyle={/S/U/W 0}]{hyperref}
\hypersetup{
    colorlinks=true,
    linkcolor=magenta,
    citecolor=cyan,
}

\usepackage{cleveref}

\usepackage{etoolbox}

\numberwithin{equation}{section}

\ifdefined\thmcolor
\declaretheoremstyle[
  shaded={bgcolor=\thmcolor}
]{plain}
\else
\fi

\ifdefined\defcolor
\declaretheoremstyle[
  headfont=\normalfont\bfseries,
  bodyfont=\normalfont,
  shaded={bgcolor=\defcolor}
]{noital}
\else
\declaretheoremstyle[
  headfont=\normalfont\bfseries,
  bodyfont=\normalfont,
]{noital}
\fi

\declaretheorem[style=plain,numberwithin=section,name=Theorem]{theorem}

\declaretheorem[style=plain,sibling=theorem,name=Conjecture]{conjecture}
\declaretheorem[style=plain,sibling=theorem,name=Claim]{claim}

\declaretheorem[style=plain,numbered=no,name=Theorem]{theorem-n}
\declaretheorem[style=plain,numbered=no,name=Proposition]{proposition-n}
\declaretheorem[style=plain,numbered=no,name=Lemma]{lemma-n}
\declaretheorem[style=plain,numbered=no,name=Corollary]{corollary-n}
\declaretheorem[style=plain,numbered=no,name=Conjecture]{conjecture-n}
\declaretheorem[style=plain,numbered=no,name=Claim]{claim-n}
\declaretheorem[style=plain,numbered=no,name=Fact]{fact-n}
\declaretheorem[style=plain,numbered=no,name=Open Problem]{openproblem-n}

\declaretheorem[style=noital,numbered=no,name=Remark]{remark-n}
\declaretheorem[style=noital,numbered=no,name=Definition]{definition-n}
\declaretheorem[style=noital,numbered=no,name=Construction]{construction-n}
\declaretheorem[style=noital,numbered=no,name=Example]{example-n}

\newcommand{\defined}{\mathrel{\coloneqq}}

\newcommand{\st}{\mathbin{\colon}}

\undef{\set}
\DeclarePairedDelimiter{\set}{\lbrace}{\rbrace}

\undef{\emptyset}
\newcommand{\emptyset}{\varnothing}

\newcommand{\union}{\mathbin{\cup}}
\newcommand{\inter}{\mathbin{\cap}}

\newcommand{\from}{\colon}

\newcommand{\divides}{\mathrel{\mid}}

\undef{\abs}
\DeclarePairedDelimiterX{\abs}[1]
  {\lvert}{\rvert}{\ifblank{#1}{\,\cdot\,}{#1}}

\undef{\norm}
\DeclarePairedDelimiterX{\norm}[1]
  {\lVert}{\rVert}{\ifblank{#1}{\,\cdot\,}{#1}}

\DeclarePairedDelimiterX{\inner}[2]
  {\langle}{\rangle}{\ifblank{#1}{\,\cdot\,}{#1},\ifblank{#2}{\,\cdot\,}{#2}}

\DeclareMathDelimiter{\given}
  {\mathbin}{symbols}{"6A}{largesymbols}{"0C}

\DeclareMathOperator{\Prob}{\mathbb{P}}
\DeclarePairedDelimiterXPP{\prob}[1]
  {\Prob}{\lparen}{\rparen}{}
  {\renewcommand{\given}{\nonscript\;\delimsize\vert\nonscript\;\mathopen{}}#1}

\DeclareMathOperator{\Expec}{\mathbb{E}}
\DeclarePairedDelimiterXPP{\expec}[1]
  {\Expec}{\lparen}{\rparen}{}
  {\renewcommand{\given}{\nonscript\;\delimsize\vert\nonscript\;\mathopen{}}#1}

\DeclareMathOperator{\Var}{Var}
\DeclarePairedDelimiterXPP{\var}[1]
  {\Var}{\lparen}{\rparen}{}
  {\renewcommand{\given}{\nonscript\;\delimsize\vert\nonscript\;\mathopen{}}#1}

\DeclareMathOperator{\Cov}{Cov}
\DeclarePairedDelimiterXPP{\cov}[2]
  {\Cov}{\lparen}{\rparen}{}{#1,#2}

\newcommand{\eps}{\varepsilon}
\newcommand{\sseq}{\subseteq}

\newcommand{\rS}{\mathrm{S}}
\newcommand{\rT}{\mathrm{T}}

\newcommand{\bfp}{\mathbf{p}}
\newcommand{\bfq}{\mathbf{q}}
